\documentclass[11pt]{amsart}
\usepackage{amssymb}
\usepackage{amscd}

\newtheorem{Thm}{Theorem}[section]
\newtheorem{Lem}[Thm]{Lemma}
\newtheorem{Cor}[Thm]{Corollary}

\newcommand{\field}{k}
\newcommand{\M}{{\rm M}}

\begin{document}

\title{Irreducible components of varieties of modules}
\author{William Crawley-Boevey}
\address{William Crawley-Boevey\newline
Department of Pure Mathematics,\newline
University of Leeds,\newline
Leeds LS2 9JT,UK}
\email{w.crawley-boevey@leeds.ac.uk}
\author{Jan Schr\"oer}
\address{Jan Schr\"oer\newline
Department of Pure Mathematics\newline
University of Leeds\newline
Leeds LS2 9JT,UK}
\email{jschroer@maths.leeds.ac.uk}


\begin{abstract}
We prove some basic results about irreducible components
of varieties of modules for an arbitrary finitely generated 
associative algebra. 
Our work generalizes results of Kac and Schofield on representations 
of quivers, but our methods are quite different, being based on deformation theory.
\end{abstract}

\maketitle

\section{Introduction}

Let $\field$ be an algebraically closed field. 
When studying the representation theory of a finitely generated 
$\field$-algebra $A$ (associative, with unit), one would often like
to classify the finite-dimensional $A$-modules. In general, however,
this is a hopeless---or meaningless---task.
For example, the 1-dimensional 
modules for a f.g.\ commutative algebra correspond to points of the 
corresponding affine variety, and what should it mean to classify the points?
A more basic problem in that case is to compute the irreducible components
of the variety. This suggests to study, for each $d\ge 1$, the irreducible 
components of the `module variety' 
${\rm mod}_A^d(\field)$ of $d$-dimensional $A$-modules.

More precisely, ${\rm mod}_A^d(\field)$ is the set of $A$-module 
structures on $k^d$, or equivalently the set of $\field$-algebra 
homomorphisms from $A$ to $\M_d(\field)$. Now a homomorphism is 
determined by its value on a set of generators of $A$, and this gives an 
embedding of ${\rm mod}_A^d(\field)$ as a closed subvariety of a suitable 
product of copies of $\M_d(\field)$. The group ${\rm GL}_d(\field)$ 
acts by conjugation on ${\rm mod}_A^d(\field)$, and the orbits 
correspond to isomorphism classes of $d$-dimensional $A$-modules. 
By abuse of language, if $C$ is a ${\rm GL}_d(\field)$-stable subset 
of ${\rm mod}_A^d(\field)$, we say that a $d$-dimensional $A$-module 
$M$ is `in' $C$ if the corresponding orbit is contained in $C$.

In addition to computing the irreducible components $C$ of 
${\rm mod}_A^d(\field)$, one would like to be able to say
something about the properties of a general module contained in $C$. 
That is, the properties which hold for all modules in some dense open
subset of $C$. For example, is the general module in $C$ indecomposable?
Can one say anything about the possible submodules of a general module in $C$?
What is the general dimension of the space of homomorphisms
or extensions between modules in two irreducible components?
In fact the dimensions of these spaces are known to be
upper semicontinuous functions (see Lemma \ref{l:uppersemihomext}), 
so given irreducible components
$C_1$ in ${\rm mod}_A^{d_1}(\field)$
and $C_2$ in ${\rm mod}_A^{d_2}(\field)$
this amounts to determining
the numbers
\[ 
\begin{split}
{\rm hom}_A(C_1,C_2) &= \min \{ {\rm dim}\, {\rm Hom}_A(M_1,M_2) \}\text{, and}
\\
{\rm ext}^1_A(C_1,C_2) &= \min \{ {\rm dim}\, {\rm Ext}^1_A(M_1,M_2) \}
\end{split}
\]
where the modules $M_i$ range over all modules in $C_i$.

The problem we consider is motivated by the fact that many
natural varieties occur as, or are related to, module varieties.
This includes the variety of complexes \cite{CS}, Lusztig's
nilpotent variety \cite{L}, \S 12, and Kleinian singularities and their
deformations \cite{CH}. The case of Lusztig's nilpotent variety is
particularly interesting, since by work of Kashiwara and
Saito \cite{KS}, its irreducible components are in 1-1 correspondence
with elements of the crystal basis of a quantum group.
Because these varieties arise as module varieties, one can
hope to use decomposition properties of modules, and
homological algebra techniques, to study the irreducible 
components. It is this theory that we initiate here.

The examples mentioned actually require a seemingly more general setup.
Fixing orthogonal idempotents $e_1,\dots,e_n\in A$ with
$e_1+\dots+e_n=1$, and a dimension vector 
$\mathbf{d}=(d_1,\dots,d_n)\in\mathbb{N}^n$,
one can consider the set ${\rm mod}_A^{\mathbf{d}}(\field)$ of
$\field$-algebra homomorphisms $A\to \M_{d_1+\dots+d_n}(\field)$,
sending each $e_i$ to the matrix whose diagonal
$d_i\times d_i$ block is the identity, and
with all other blocks zero.
The natural group which acts is 
the product of the general linear groups of size $d_i$.
Taking $n=1$ one recovers ${\rm mod}_A^d(\field)$.
Taking $A=kQ$, the path algebra of a quiver $Q$, and letting the $e_i$
be the `trivial paths', one obtains the space ${\rm Rep}(Q,\mathbf{d})$
of representations of $Q$.
All of our results below hold in this more general setup, either
by adapting the proofs, or by using Bongartz's observation \cite{B1} that 
these varieties are related by a fibre bundle construction.

In fact ${\rm Rep}(Q,\mathbf{d})$ is a model for some of our work.
There is no problem describing the irreducible 
components in this case, as ${\rm Rep}(Q,\mathbf{d})$ 
is a vector space, hence irreducible. But the problem of 
determining when the general element is an indecomposable 
representation, and of computing the general dimension
of homomorphism and extension spaces is nontrivial.
A nice theory has, however, been created by 
Kac \cite{K} and by Schofield \cite{S}.
Our more general setup differs from this one,
however, in an important way: it is possible to
relate certain Ext spaces to tangent spaces
in ${\rm Rep}(Q,\mathbf{d})$, but in 
${\rm mod}_A^d(\field)$ this is no longer the case.
This is because Voigt's Lemma \cite{Ga} really 
involves a certain scheme ${\rm mod}_A^d$
which need not be reduced, not even generically reduced, 
so its tangent spaces are not the same as for ${\rm mod}_A^d(\field)$. 

Given subsets $C_i \subseteq {\rm mod}_A^{d_i}(\field)$ 
which are ${\rm GL}_{d_i}(\field)$-stable ($1\le i\le t$), 
we consider all modules of dimension $d = d_1+\dots+d_t$
which are of the form $M_1\oplus \dots \oplus M_t$ with 
the $M_i$ in $C_i$, and we denote by 
$C_1\oplus\dots\oplus C_t$ the corresponding
${\rm GL}_{d}(\field)$-stable subset of ${\rm mod}_A^{d}(\field)$.
This is the image of the map
\begin{equation}
\label{e:mapeta}
{\rm GL}_d(\field) \times C_1 \times \cdots \times C_t
\longrightarrow {\rm mod}_A^d(\field) 
\end{equation}
sending a tuple $(g,x_1,\dots,x_t)$ to the conjugation by $g$
of the module structure in ${\rm mod}_A^d(\field)$ which has
the $x_i$ as diagonal blocks.
It follows that if the $C_i$ are irreducible locally closed subsets, then the
closure $\overline{C_1\oplus\dots\oplus C_t}$ is irreducible.

Our first result is an analogue of the Krull-Remak-Schmidt
Theorem, and indeed follows from it quite easily.
Most of the theorem is already known to the experts, but the only reference 
seems to be the preprint \cite{JAP}, which is not published elsewhere.
For convenience we give a full proof.

\begin{Thm}\label{th:KS}
If $C$ is an irreducible component in ${\rm mod}_A^{d}(\field)$,
then 
\begin{equation}
\label{e:candec}
C = \overline{C_1\oplus\dots\oplus C_t}
\end{equation}
for some irreducible components $C_i$ of module varieties 
${\rm mod}_A^{d_i}(\field)$, 
with the property that the general module in each $C_i$ is 
indecomposable. Moreover $C_1,\dots,C_t$
are uniquely determined by this, up to reordering.
\end{Thm}

One cannot use direct sums quite as freely as this suggests,
however, as the closure of a direct sum of irreducible
components is not in general an irreducible component.
Our main result, proved using deformation theory, 
is the determination of when this happens.

\begin{Thm}\label{th:sum}
If $C_i \subseteq {\rm mod}_A^{d_i}(\field)$ are irreducible
components ($1\le i\le t$), and $d = d_1+\dots+d_t$,
then $\overline{C_1\oplus\dots\oplus C_t}$
is an irreducible component of ${\rm mod}_A^d(\field)$
if and only if ${\rm ext}_A^1(C_i,C_j) = 0$ for all $i \not= j$.
\end{Thm}

These two theorems give an analogue of
Kac's `canonical decomposition' for representations of
quivers \cite{K}, Proposition 3. Thus we call (\ref{e:candec}) the 
\emph{canonical decomposition} of $C$. It reduces the problem
of computing all irreducible components of module varieties 
to the case of components
whose general element is indecomposable.  
By studying the fibres of the map
(\ref{e:mapeta}), it is easy to see that 
\[
\dim C = \sum_{i=1}^t \dim C_i + \sum_{i \not= j} 
\left(d_i d_j - {\rm hom}_A(C_i,C_j)\right),
\]
and also, if $D$ is any irreducible component of a module variety, then
\[
{\rm hom}_A(C,D) = \sum_{i=1}^t {\rm hom}_A(C_i,D),
\]
and similarly for ${\rm ext}^1_A$, or with $C$ and $D$ exchanged.

Instead of taking direct sums of the modules in two irreducible
components, one can take extensions. 
Let $d = d_1+d_2$, let
\[
G = {\rm GL}_{d_1}(\field)\times{\rm GL}_{d_2}(\field),
\]
and let $S$ be a $G$-stable subset of 
${\rm mod}_A^{d_1}(\field)\times {\rm mod}_A^{d_2}(\field)$.
Again, by abuse of language we say that a pair
of modules $(M_1,M_2)$ is `in' $S$ if the product of the
corresponding orbits is a subset of $S$.
We denote by ${\mathcal E}(S)$
the ${\rm GL}_{d}(\field)$-stable subset of ${\rm mod}_A^{d}(\field)$ 
corresponding to all modules $M$ which belong to
a short exact sequence
\begin{equation}
\label{e:ses}
0\longrightarrow M_2 \longrightarrow M\longrightarrow M_1\longrightarrow 0
\end{equation}
with $(M_1,M_2)$ in $S$. We have the following result.

\begin{Thm}
\label{th:multi}
Let $S$ be a $G$-stable subset of 
${\rm mod}_A^{d_1}(\field)\times {\rm mod}_A^{d_2}(\field)$.
\begin{itemize}
\item[(i)]If $S$ is a closed subset, then ${\mathcal E}(S)$ 
is a closed subset of ${\rm mod}_A^{d}(\field)$.

\item[(ii)]If $S$ is an irreducible locally closed subset, and
${\rm Ext}^1_A(M_1,M_2)$ has constant dimension for all $(M_1,M_2)$ 
in $S$, then the closure of ${\mathcal E}(S)$ in ${\rm mod}_A^{d}(\field)$
is irreducible.

\item[(iii)]If $S$ is an open subset and ${\rm Ext}^1_A(M_2,M_1)=0$ 
for all $(M_1,M_2)$ in $S$, then ${\mathcal E}(S)$ is an open subset
of ${\rm mod}_A^{d}(\field)$.
\end{itemize}
\end{Thm}

This has the following application to irreducible components.

\begin{Cor}
\label{co:min}
If $C_1$ and $C_2$ are irreducible components of 
${\rm mod}_A^{d_1}(\field)$ and ${\rm mod}_A^{d_2}(\field)$ respectively,
and if ${\rm ext}_A^1(C_2,C_1) = 0$, then
${\mathcal E}(C_1\times C_2)$ contains an irreducible 
component of ${\rm mod}_A^{d}(\field)$.
It may be realized as $\overline{{\mathcal E}(S)}$, 
where $S$ is the $G$-stable subset 
of $C_1\times C_2$
corresponding to the pairs of modules $(M_1,M_2)$
with $\dim {\rm Ext}_A^1(M_2,M_1) = 0$ and
$\dim {\rm Ext}_A^1(M_1,M_2) = {\rm ext}_A^1(C_1,C_2)$.
\end{Cor}

In particular, if ${\rm mod}_A^{d}(\field)$ is known to be irreducible
then, under the hypotheses of the theorem, associated to 
every $d$-dimensional $A$-module $M$ there
must be a short exact sequence (\ref{e:ses}) with $M_i$ in $C_i$. 
Our results thus generalize the implications 
(i)$\Leftrightarrow$(ii)$\Leftarrow$(iii) 
of \cite{S}, Theorem 3.3. The remaining implication 
is clearly false in general, even for 2-dimensional modules
for the algebra of dual numbers.

One source of irreducible components is provided by modules $M$ without
self-extensions, i.e.\ ${\rm Ext}^1_A (M,M)=0$.
This condition ensures that the corresponding orbit ${\mathcal O}(M)$ is open
in ${\rm mod}_A^d(\field)$, see \cite{Ga}, and so, since it is irreducible, 
its closure $\overline{{\mathcal O}(M)}$ is an irreducible component. 

It is perhaps worthwhile to record the following consequence
of Theorem \ref{th:multi}(iii). If $M_1,\dots,M_t$ are $A$-modules, 
we denote by
${\mathcal F}(M_1,\dots,M_t)$ the set of $A$-modules $N$ which have an
increasing sequence 
$0 = N_0 \subseteq N_1 \subseteq \cdots \subseteq N_s = N$
of submodules such that each quotient $N_j/N_{j-1}$ is isomorphic to one of the $M_i$.
We denote by ${\mathcal F}_d(M_1, \cdots, M_t)$ the 
${\rm GL}_d(\field)$-stable subset of ${\rm mod}_A^d(\field)$
corresponding to the $d$-dimensional modules in 
${\mathcal F}(M_1,\dots,M_t)$.

\begin{Cor}
\label{co:filt}
If $M_1, \cdots, M_t$ are modules with ${\rm Ext}_A^1(M_i,M_j) = 0$
for all $i \leq j$, then ${\mathcal F}_d(M_1, \cdots, M_t)$
is open in ${\rm mod}_A^d(\field)$ for all $d$.
\end{Cor}

In this paper, unless otherwise stated, `modules' are right modules.
Although we often write maps on the left hand side,
we compose them as if they were on the right. 
Thus the composition of a map $\theta$ 
followed by a map $\phi$ is denoted $\theta\phi$.
By a variety we mean a locally closed subset of affine or projective space
over $\field$.
We prove the results in numerical order, except that 
Theorem \ref{th:sum} is proved last.
At the end we give various examples.

\section{Proof of Theorem \ref{th:KS}}

If $C$ is an irreducible component of ${\rm mod}_A^d(\field)$, 
we denote by $C^{\rm ind}$ the subset of $C$ given by the
indecomposable modules. This is a constructible subset of
${\rm mod}_A^d(\field)$, and clearly every indecomposable 
$d$-dimensional module is in $C^{\rm ind}$ for some
irreducible component $C$.

Every $d$-dimensional module $M$ is isomorphic to a
direct sum of indecomposables, and if the summands have 
dimensions $d_1,\dots,d_t$ then $M$ is in a subset
\[
S = C_1^{\rm ind} \oplus \dots \oplus C_t^{\rm ind}
\]
for some irreducible components 
$C_i\subseteq {\rm mod}_A^{d_i}(\field)$.
Since the sets $S$ which arise this way are constructible
in ${\rm mod}_A^d(\field)$, and their union is the whole space, 
if $C$ is an irreducible component of ${\rm mod}_A^d(\field)$,
at least one of the sets $S$ must contain a dense open subset of $C$.
We concentrate on this set.

As observed in the introduction, $\overline{C_1\oplus\dots\oplus C_t}$ 
is an irreducible closed subset of ${\rm mod}_A^d(\field)$.
Since it contains $C$, which is an irreducible component, it follows that
\[
C = \overline{C_1\oplus\dots\oplus C_t}.
\]
Now suppose that
\[
(g,x_1,\dots,x_t) \in {\rm GL}_d(\field) \times C_1 \times \cdots \times C_t
\]
is an element which is sent under the map (\ref{e:mapeta}) to
an element of $S$.
By the Krull-Remak-Schmidt Theorem, any two decompositions of a module
into a direct sum of indecomposables have the same number of
terms. It follows that $x_i\in C_i^{\rm ind}$ for all $i$. Thus the
inverse image of $S$ is
${\rm GL}_d(\field) \times C_1^{\rm ind} \times \cdots \times C_t^{\rm ind}$.
Since $C_1\oplus\dots\oplus C_t$ contains a dense open subset of $C$,
${\rm GL}_d(\field) \times C_1^{\rm ind} \times \cdots \times C_t^{\rm ind}$ 
contains a dense open subset of 
${\rm GL}_d(\field) \times C_1 \times \cdots \times C_t$.
It follows that each $C_i^{\rm ind}$ contains a dense open subset of $C_i$.
This completes the proof of existence.

Now uniqueness. Suppose that $C$ is an irreducible component with
\[
C = \overline{C_1\oplus \dots \oplus C_t} = \overline{D_1\oplus \dots \oplus D_s}
\]
for irreducible components $C_i$ and $D_j$ with the general element indecomposable.

We denote by $C_i^0$ the set of elements of $C_i$ which are not contained in
any other irreducible component of ${\rm mod}_A^{d_i}(\field)$. 
Clearly $C_i^0$ is a dense open subset of $C_i$. Now the general element of $C_i$ is
indecomposable, so $C_i^{0,\rm ind} = C_i^0 \cap C_i^{\rm ind}$ contains
a dense open subset of $C_i$. 
Thus 
\[
{\rm GL}_d(\field) \times C_1^{0,\rm ind} \times \cdots \times C_t^{0,\rm ind}
\]
contains a dense open subset of
\[
{\rm GL}_d(\field) \times C_1 \times \cdots \times C_t.
\]
Now if $X\to Y$ is a dominant morphism between irreducible varieties
then the image of a dense open subset of $X$ must contain a dense 
open subset of $Y$. It follows that 
$C_1^{0,\rm ind} \oplus \cdots \oplus C_t^{0,\rm ind}$
contains a dense open subset of $C$.
Repeating for the $D_j$, and using the fact that two dense open
subsets of $C$ must intersect, we obtain an element in
\[
(C_1^{0,\rm ind} \oplus \cdots \oplus C_t^{0,\rm ind})
\:
\cap
\:
(D_1^{0,\rm ind} \oplus \cdots \oplus D_s^{0,\rm ind}).
\]
Thus there is an isomorphism of modules
\[
M_1\oplus\dots\oplus M_t \cong N_1\oplus\dots\oplus N_s
\]
with the $M_i$ in $C_i^{0,\rm ind}$ and the $N_j$ in $D_j^{0,\rm ind}$.
By the Krull-Remak-Schmidt Theorem, we have $s=t$, and after permuting
the $N_j$ we have $M_i\cong N_i$ for all $i$. It follows that $C_i = D_i$.

\section{Proof of Theorem \ref{th:multi}{\rm (i)}}

Let $S$ be a closed subset of
${\rm mod}_A^{d_1}(\field)\times {\rm mod}_A^{d_2}(\field)$, 
stable for the action of $G = {\rm GL}_{d_1}(\field)\times {\rm GL}_{d_2}(\field)$,
and let $d = d_1 + d_2$.
By $V^{\rm ses}(d_1,d_2)$ we denote the locally closed subset of pairs
$(\theta,\phi)$ in ${\rm Hom}_\field(\field^{d_2},\field^d)
\times {\rm Hom}_\field(\field^{d},\field^{d_1})$ such that
\[ 
0 \longrightarrow \field^{d_2} \stackrel{\theta}{\longrightarrow} 
\field^d \stackrel{\phi}{\longrightarrow} \field^{d_1} 
\longrightarrow 0 
\]
is a short exact sequence.

Let $V_A^{\rm ses}$ be the set of elements
$(m,(\theta,\phi))$ in
${\rm mod}_A^d(\field) \times V^{\rm ses}(d_1,d_2)$ such that
the image of $\theta$ is a submodule of $\field^d$, with
the module structure given by $m$. 
Recall that $m$ is an algebra
homomorphism $A \to \M_d(\field)$.
Identifying elements of $\M_d(\field)$ with endomorphisms
of $\field^d$, the condition is that 
$\theta m(a) \phi = 0$ for all $a$ in $A$.
Thus $V_A^{\rm ses}$ is a closed subset of 
${\rm mod}_A^d(\field) \times V^{\rm ses}(d_1,d_2)$.

Let
\[ 
\mu: V_A^{\rm ses} \longrightarrow {\rm mod}_A^{d_1}(\field)\times {\rm mod}_A^{d_2}(\field)
\]
be the map sending an element $(m,(\theta,\phi))$ to the pair 
$(m_1,m_2)$ where 
$m_1$ and $m_2$ are the unique module
structures on $\field^{d_1}$ and $\field^{d_2}$ for which $\theta$
and $\phi$ are module homomorphisms, where $\field^d$ is
considered as an $A$-module using the structure $m$.

It is not obvious from this, but it can be checked, that $\mu$
is a morphism of varieties. One uses the covering
of $V^{\rm ses}(d_1,d_2)$ by affine open subsets
consisting of the pairs $(\theta,\phi)$ where the image of $\theta$ is 
complementary to a given subspace of $\field^d$.
See \cite{R} for a similar situation.

Now $\mu^{-1}(S)$ is a closed subset of $V_A^{\rm ses}$, and
hence a closed subset of 
${\rm mod}_A^d(\field) \times V^{\rm ses}(d_1,d_2)$.

By ${\rm Gr}(d_2 \hookrightarrow d)$ we denote
the Grassmannian of $d_2$-dimensional subspaces of
$\field^d$.
Let 
\[ 
\pi: V^{\rm ses}(d_1,d_2) \longrightarrow {\rm Gr}(d_2 \hookrightarrow d) 
\] 
be the map which sends a pair 
$(\theta,\phi)$ to the image of $\theta$.
Clearly $\pi$ is a morphism of varieties, and it is a principal 
$G$-bundle.
Since $S$ is $G$-stable we know that $\mu^{-1}(S)$ is a union of fibres 
of the morphism 
\[ 
(1,\pi): {\rm mod}_A^d(\field) \times V^{\rm ses}(d_1,d_2)
\longrightarrow {\rm mod}_A^d(\field) \times 
{\rm Gr}(d_2 \hookrightarrow d). 
\]
Now the principal bundle property implies that the image
$(1,\pi)(\mu^{-1}(S))$ must be closed
in ${\rm mod}_A^d(\field) \times {\rm Gr}(d_2 \hookrightarrow d)$.

Let 
\[ 
\pi_1:  {\rm mod}_A^d(\field) \times 
{\rm Gr}(d_2 \hookrightarrow d) \longrightarrow {\rm mod}_A^d(\field) 
\]
be the first projection.
Since ${\rm Gr}(d_2 \hookrightarrow d)$ is a projective variety,
the image of  $(1,\pi)(\mu^{-1}(S))$ under
$\pi_1$ is closed.
Now this image is clearly ${\mathcal E}(S)$,
and the result is proved.

\section{Derivations}

Recall that if $M$ is an $A$-$A$-bimodule, then a \emph{derivation} 
$d:A\to M$ is a
linear map with $d(ab) = ad(b)+d(a)b$ for all $a,b\in A$. We denote by
${\rm Der}(A,M)$ the vector space of all derivations $A\to M$. 
A derivation $d$ is
\emph{inner} if there is some $m\in M$ with $d(a) = am-ma$ for all $a\in A$.
We denote the subspace of inner derivations by ${\rm Der}^0(A,M)$.

\begin{Lem}
Let $a_1,\dots,a_N$ be $\field$-algebra generators of $A$.
If $M$ is an $A$-$A$-bimodule, then
${\rm Der}(A,M)$ is isomorphic to the subspace of $M^N$
consisting of the tuples $(d_1,\dots,d_N)$ with the property that
\begin{equation}
\label{e:dercond}
\sum_{i_1,\dots,i_r} \lambda_{i_1,\dots,i_r} 
\left(\sum_{j=1}^r a_{i_1} \dots a_{i_{j-1}} d_{i_j} a_{i_{j+1}} \dots a_{i_r}\right) = 0
\end{equation}
for all noncommutative polynomials
\[
f(x_1,\dots,x_N) = \sum_{i_1,\dots,i_r} \lambda_{i_1,\dots,i_r} x_{i_1} \dots x_{i_r}
\]
over $\field$, with the property that $f(a_1,\dots,a_N)=0$.
\end{Lem}

\begin{proof}
A derivation $d$ is determined by its values on the generators $a_i$,
and if $d(a_i) = d_i$, then
\begin{equation}
\label{e:derprodform}
d(a_{i_1} \dots a_{i_r}) = \sum_{j=1}^r a_{i_1} \dots a_{i_{j-1}} d_{i_j} a_{i_{j+1}} \dots a_{i_r}.
\end{equation}
Thus if $f(a_1,\dots,a_N)=0$, then $d(f(a_1,\dots,a_N))=0$, which
implies (\ref{e:dercond}). 
Conversely, given $d_1,\dots,d_N$, one can define a function
$d:A\to M$ using (\ref{e:derprodform}). The condition (\ref{e:dercond})
ensures that this is well-defined. The construction of $d$ ensures that
it is a derivation.
\end{proof}

Recall that if $X$ is a variety, then a function $f:X\to \mathbb{Z}$ is \emph{upper semicontinuous}
if $\{ x\in X | f(x)\le n\}$ is open in $X$ for all $n\in\mathbb{Z}$.
If $V$ is a vector space, then a \emph{cone} in $V$ is a subset which contains 0 and
is closed under multiplication by elements of the field $\field$. 
The following result is well-known.

\begin{Lem}
\label{l:uppersemi}
Let $X$ be a variety, $V$ is a vector space and $F$ a closed subset
of $X\times V$. If, for all $x\in X$ the set $F_x = \{ v\in V | (x,v)\in F\}$
is a subspace of $V$, or more generally a cone in $V$, 
then the function $X\to\mathbb{Z}$ sending $x$ to
$\dim F_x$ is upper semicontinuous.
\end{Lem}

If $M_1$ and $M_2$ are $A$-modules, then ${\rm Hom}_\field(M_1,M_2)$
is naturally an $A$-$A$-bimodule, and it is well-known that
\begin{equation}
\label{e:extder}
{\rm Ext}^1_A(M_1,M_2) \cong 
{\rm Der}(A,{\rm Hom}_\field(M_1,M_2)) / 
{\rm Der}^0(A,{\rm Hom}_\field(M_1,M_2)).
\end{equation}
From the construction we also have
\[
{\rm Der}^0(A,{\rm Hom}_\field(M_1,M_2)) \cong {\rm Hom}_\field(M_1,M_2) / 
{\rm Hom}_A (M_1,M_2),
\]
and hence
\begin{multline}
\label{e:extderform}
\dim {\rm Ext}^1_A(M_1,M_2) = \dim {\rm Der}(A,{\rm Hom}_\field(M_1,M_2)) +\\
+ \dim {\rm Hom}_A (M_1,M_2) - \dim M_1 \dim M_2.
\end{multline}

\begin{Lem}
\label{l:uppersemihomext}
The functions
\[
{\rm mod}_A^{d_1}(\field)\times {\rm mod}_A^{d_2}(\field)\longrightarrow 
\mathbb{Z}
\]
sending a pair $(m_1,m_2)$ to the dimensions of the 
spaces ${\rm Hom}_A(M_1,M_2)$,
${\rm Der}(A,{\rm Hom}_\field (M_1,M_2))$ 
and ${\rm Ext}^1_A(M_1,M_2)$ are upper semicontinuous
(where $M_i$ is the module corresponding to $m_i$).
\end{Lem}

For ${\rm Hom}$ spaces this is well-known. For ${\rm Ext}^1$ it
is known, but we could not find a convenient reference
which applies for all finitely generated $\field$-algebras $A$,
and the usual method, using projective resolutions, does not
immediately adapt. 

\begin{proof}
For ${\rm Hom}_A(M_1,M_2)$ one applies Lemma \ref{l:uppersemi} 
with the set $F_1$ of 
\[
((m_1,m_2),\theta) \in ({\rm mod}_A^{d_1}(\field)\times {\rm mod}_A^{d_2}(\field)) 
\times {\rm Hom}_\field (\field^{d_1},\field^{d_2})
\]
for which $\theta$ is an $A$-module homomorphism, for the module
structures given by $m_1$ and $m_2$. 

For ${\rm Der}(A,{\rm Hom}_\field (M_1,M_2))$ 
one uses the set $F_2$ of
\[
((m_1,m_2),(d_1,\dots,d_N)) \in ({\rm mod}_A^{d_1}(\field)\times {\rm mod}_A^{d_2}(\field)) 
\times {\rm Hom}_\field (\field^{d_1},\field^{d_2})^N
\]
for which $(d_1,\dots,d_N)$ satisfies the condition (\ref{e:dercond}) using
the module structures $m_1$ and $m_2$ to turn 
${\rm Hom}_\field (\field^{d_1},\field^{d_2})$ into an $A$-$A$-bimodule.
This is clearly a closed subset.

Finally, a sum of upper semicontinuous functions is upper semicontinuous, so
the assertion for ${\rm Ext}^1_A(M_1,M_2)$ follows from
equation (\ref{e:extderform}).
\end{proof}

\begin{Lem}
\label{l:homderconst}
If $X$ is an irreducible locally closed subset of 
${\rm mod}_A^{d_1}(\field)\times {\rm mod}_A^{d_2}(\field)$
on which the dimension ${\rm Ext}^1_A(M_1,M_2)$
is constant, then so are the dimensions of
${\rm Hom}_A(M_1,M_2)$ and ${\rm Der}(A,{\rm Hom}_\field (M_1,M_2))$.
\end{Lem}

\begin{proof}
The subsets of $X$ on which the dimensions of
${\rm Hom}_A(M_1,M_2)$ and ${\rm Der}(A,{\rm Hom}_\field (M_1,M_2))$
take their minimum values are nonempty open subsets of $X$,
and so they must meet. Away from this intersection ${\rm Ext}^1_A(M_1,M_2)$
must have non-minimal dimension.
\end{proof}

\section{Proof of Theorem \ref{th:multi}{\rm (ii)}}

Suppose that $S$ is an irreducible locally closed $G$-stable subset of
$X = {\rm mod}_A^{d_1}(\field)\times {\rm mod}_A^{d_2}(\field)$,
and suppose that ${\rm Ext}^1_A(M_1,M_2)$ has dimension $e$
for all $(M_1,M_2)$ in $S$. Let $d=d_1+d_2$.

Let $Z$ be the closed subset of ${\rm mod}_A^d(\field)$ consisting of the
homomorphisms $A\to \M_d(\field)$ which take an upper
triangular block form for the decomposition $d=d_1+d_2$ with
the lower $d_2\times d_1$ block zero.
A map $m:A\to \M_d(\field)$ taking this block form is a homomorphism
if and only if the diagonal blocks are homomorphisms $m_i:A\to \M_{d_i}(\field)$
and the upper triangular block 
$A \to \M_{d_1\times d_2}(\field)$
is a derivation for the $A$-$A$-bimodule structure
on $\M_{d_1\times d_2}(\field) \cong {\rm Hom}_\field(\field^{d_1},\field^{d_2})$
given by the module structures $m_1$ and $m_2$.
Thus $Z$ is isomorphic to the closed subset $F_2$ used in the
proof of Lemma \ref{l:uppersemihomext}.

By Lemma \ref{l:homderconst}, the space 
${\rm Der}(A,{\rm Hom}_\field (M_1,M_2))$
has constant dimension on $S$. It follows that
the projection 
$\pi:(S \times {\rm Hom}_\field(\field^{d_1},\field^{d_2})^N) \cap F_2 \to S$ is a 
sub-bundle of the trivial bundle 
$S \times {\rm Hom}_\field (\field^{d_1},\field^{d_2})^N \to S$.
In particular it is a vector bundle in its own right. 
Thus, since $S$ is irreducible, so is $\pi^{-1}(S)$.

Now the composition 
$\pi^{-1}(S)\hookrightarrow F_2\cong Z\hookrightarrow {\rm mod}_A^d(\field)$
and the conjugation action of ${\rm GL}_d(\field)$ combine to
give a map
\[
{\rm GL}_d(\field)\times \pi^{-1}(S) \longrightarrow {\rm mod}_A^d(\field).
\]
Now the left hand side is irreducible, so its image,
which is ${\mathcal E}(S)$, has irreducible closure.

\section{Deformations}\label{deformations}

Let $\field[[T]]$ be the power series algebra in one variable over $\field$.
Given a homomorphism $m:A\to \M_d(k)$, a {\it deformation} of $m$ is a 
homomorphism $\mathbf{m}:A\to \M_d(\field[[T]])$ whose composition with the natural 
homomorphism $\M_d(\field[[T]])\to \M_d(k)$ is equal to $m$. Equivalently, 
for each $a\in A$, the $(i,j)$ component of $m(a)$ is the constant term
of the $(i,j)$ component of $\mathbf{m}(a)$.

The next lemma follows the philosophy of Gerstenhaber's deformation
theory, see for example \cite{Ge} or \cite{GS}. Let $d_1$ and $d_2$
be positive integers and set $d = d_1+d_2$.

\begin{Lem}
\label{l:deform}
Let $m:A\to \M_d(k)$ be a homomorphism which takes the
upper triangular block form
\[
m(a) = \begin{pmatrix} m_{11}(a) & m_{12}(a) \\ 0 & m_{22}(a) \end{pmatrix},
\]
where $m_{ij}(a)$ is a matrix of size $d_i\times d_j$, and let $M_i$ 
be the $A$-module of dimension $d_i$ corresponding to $m_{ii}$.
Suppose ${\rm Ext}^1_A(M_2,M_1)=0$.
Then for any deformation of $m$ there exists an element $\mathbf{g}$
with block form
\[
\mathbf{g} =  \begin{pmatrix} 1 & 0 \\ \mathbf{g}_{21} & 1\end{pmatrix} \in 
{\rm GL}_d(\field[[T]]),
\]
where $\mathbf{g}_{21}$ is a $d_2\times d_1$ matrix with entries in 
$T\field[[T]]$, such that the conjugate of the deformation by $\mathbf{g}$
is a deformation
$\mathbf{m} : A\to \M_d(\field[[T]])$,
where $\mathbf{m}$ has the 
upper triangular block form
\[
\mathbf{m}(a) = \begin{pmatrix} \mathbf{m}_{11}(a) & \mathbf{m}_{12}(a) \\ 
0 & \mathbf{m}_{22}(a) \end{pmatrix}.
\]
\end{Lem}

\begin{proof}
We show first that if $\mathbf{m}$ is a deformation of $m$ and
for some $n\ge 1$ the lower triangular 
block of $\mathbf{m}(a)$ has coefficients in $T^n \field[[T]]$ for all 
$a\in A$, 
then there is a matrix 
\begin{equation}
\label{e:gn}
\mathbf{g}_n = \begin{pmatrix} 1 & 0 \\ T^n \theta_n & 1 \end{pmatrix} \in 
{\rm GL}_d (\field[[T]])
\end{equation}
where $\theta_n \in \M_{d_2\times d_1}(\field)$, 
such that the lower triangular block of 
$\mathbf{g}_n \mathbf{m}(a) \mathbf{g}_n^{-1}$ has coefficients in $T^{n+1} \field[[T]]$
for all $a\in A$.

Write 
\[
\mathbf{m}(a) = \sum_{r=0}^\infty c^r (a) T^r
\]
for some maps $c^r: A \to \M_d(\field)$, and write each $c^r$ in block form
\[ 
c^r(a) = \begin{pmatrix} c^r_{11}(a) & c^r_{12}(a) \\ c^r_{21}(a) & c^r_{22}(a) 
\end{pmatrix}. 
\]
Thus
\[
\mathbf{m}(a) = \begin{pmatrix} \mathbf{m}_{11}(a) & \mathbf{m}_{12}(a) \\ 
\mathbf{m}_{21}(a) & \mathbf{m}_{22}(a) \end{pmatrix}
\]
where $\mathbf{m}_{ij}(a) = \sum_{r=0}^\infty c^r_{ij}(a) T^r$.

By assumption the constant term 
satisfies $c^0_{11} = m_{11}$, $c^0_{22} = m_{22}$, $c^0_{12} = m_{12}$ and
$c^0_{21} = 0$, and we are assuming that $c^r_{21} = 0$ for $0\le r\le n-1$.
The fact that $\mathbf{m}$ is an algebra homomorphism implies that
\[
c^r(ab) = \sum_{\stackrel{s+t=r}{s,t \geq 0}} c^s(a)c^t(b). 
\]
for all $a,b\in A$. Compare for example \cite{GS}, p.24.
Now since $c^r_{21}=0$ for all $0\le r\le n-1$, we have
\[
c^n_{21}(ab) = c^0_{22}(a)c^n_{21}(b)  + c^n_{21}(a) c^0_{11}(b) 
= m_{22} (a)c^n_{21}(b)  + c^n_{21}(a) m_{11} (b).
\]
Identifying $\M_{d_2\times d_1}(\field)$ with 
${\rm Hom}_\field(M_2,M_1)$, and considering this
as an $A$-$A$-bimodule using the actions of $A$
on $M_1$ and $M_2$, it follows that $c^n_{21}$ defines a derivation 
$A\to{\rm Hom}_\field(M_2,M_1)$.
Now by (\ref{e:extder}) and our assumption that ${\rm Ext}^1_A(M_2,M_1)=0$, 
every derivation is inner. Therefore, there is a matrix 
$\theta_n \in \M_{d_2\times d_1}(\field)$ with
$c^n_{21}(a) = m_{22}(a) \theta_n - \theta_n m_{11}(a)$
for all $a\in A$.
Let $\mathbf{g}_n$ be the matrix given by (\ref{e:gn}). Then
\[
\mathbf{g}_n^{-1} = \begin{pmatrix} 1 & 0 \\ - T^n \theta_n & 1 \end{pmatrix}
\]
and so the lower triangular $d_2\times d_1$ block of $\mathbf{g}_n \mathbf{m}(a) \mathbf{g}_n^{-1}$
is equal to
\[
\mathbf{m}_{21}(a) + T^n\theta_n \mathbf{m}_{11}(a) - 
\mathbf{m}_{22}(a)T^n\theta_n - T^{2n}\theta_n \mathbf{m}_{12}(a)\theta_n.
\]
This has entries in $T^{n+1} \field[[T]]$, so $\mathbf{g}_n$ has the required property.

We now prove the lemma. Let $\mathbf{m}^0$ be a deformation of $m$.
For each $a\in A$, the lower triangular block of $\mathbf{m}^0(a)$ has
entries in $T \field[[T]]$, so by the construction above there is
a matrix $\mathbf{g}_1$ with the property that if $\mathbf{m}^1$ is defined
by $\mathbf{m}^1(a) = \mathbf{g}_1 \mathbf{m}^0(a) \mathbf{g}_1^{-1}$
then the lower triangular block of $\mathbf{m}^1(a)$ has entries
in $T^2 \field [[T]]$. Thus there is a $\mathbf{g}_2$ such that
the lower triangular block of
$\mathbf{m}^2(a) = \mathbf{g}_2 \mathbf{m}^1(a) \mathbf{g}_2^{-1}$
has entries in $T^3 \field[[T]]$. Repeating in this way
gives deformations $\mathbf{m}^1,\mathbf{m}^2,\dots$
of $m$ and matrices $\mathbf{g}_1,\mathbf{g}_2,\dots$ of the form (\ref{e:gn}).
Define
\[
\mathbf{g} = \begin{pmatrix} 1 & 0 \\ \sum_{n=1}^\infty T^n \theta_n & 1 \end{pmatrix} \in {\rm GL}_d (\field[[T]]).
\]
and define $\mathbf{m}$ by $\mathbf{m}(a) = \mathbf{g} \mathbf{m}^0(a) \mathbf{g}^{-1}$.

We have
\[
\mathbf{g} = \mathbf{g}_n \dots \mathbf{g}_1 + O(T^{n+1})
\]
where the notation $O(T^{n+1})$ means we ignore powers $T^{n+1}$ and above. Thus
\[
\mathbf{m}(a) = \mathbf{g}_n \dots \mathbf{g}_1\mathbf{m}^0(a) \mathbf{g}_1^{-1}\dots\mathbf{g}_n^{-1} + O(T^{n+1})
= \mathbf{m}^n(a) + O(T^{n+1}).
\]
Thus, for all $n$, the lower triangular block of $\mathbf{m}(a)$ has entries in $T^{n+1}\field[[T]]$.
Thus it is zero.
\end{proof}

In another language, a \emph{deformation} of a $d$-dimensional
$A$-module $M$ is a $\field[[T]]$-$A$-bimodule $\mathbf{M}$ which
is free of rank $d$ over $\field[[T]]$ and with 
\[
\field \otimes_{\field[[T]]} \mathbf{M}\cong M.
\]
The lemma then says that
if $M$ belongs to an exact sequence
\[
0\longrightarrow M_2\longrightarrow M\longrightarrow M_1\longrightarrow 0
\]
with ${\rm Ext}^1_A(M_2,M_1)=0$, then any deformation $\mathbf{M}$ of $M$ belongs to an
exact sequence
\[
0\longrightarrow \mathbf{M}_2\longrightarrow \mathbf{M}\longrightarrow \mathbf{M}_1\longrightarrow 0
\]
with $\mathbf{M}_i$ a deformation of $M_i$.

\section{A valuative criterion}

We use the following variation on the valuative criterion
for flatness \cite{Gr}, \S 11.8, to transfer results about deformations
to module varieties. Here $\field[[T]]$ is
the power series algebra in one variable, and  
$\pi^*:{\rm Spec}(\field) \to {\rm Spec}(\field[[T]])$
is the unique morphism of schemes over $\field$.

\begin{Lem}
\label{l:schtrick}
Let $f:X\to Y$ be a morphism of schemes over $k$, and assume that
$Y$ is of finite type and quasiprojective over $k$. 
Suppose that for every commutative square of morphisms of schemes over
$\field$,
\[
\begin{CD}
{\rm Spec}(\field) @>p>> X  \\
@V\pi^*VV @VfVV \\
{\rm Spec}(\field[[T]]) @>h>> Y
\end{CD}
\]
there exists a morphism 
\[
h': {\rm Spec}(\field[[T]]) \longrightarrow X
\] 
such that $\pi^* h' = p$ and $h'f = h$, that is, the two resulting 
triangles commute.
Then the image of $f$ is open in $Y$.
\end{Lem}

\begin{proof}
Recall that the points of a scheme are in 1-1 correspondence
with irreducible closed subsets, with the point $z$
corresponding to the closure $\overline{\{z\}}$.
Since $Y$ has finite type over $\field$, any nonempty 
constructible subset of $Y$ must contain a closed point. 

Since ${\rm Im}(f)$ is a constructible subset of $Y$,
to show it is open it suffices to show that ${\rm Im}(f)$ is closed 
under generization, that is, if $z\in Y$, $y\in{\rm Im}(f)$ and
$y\in\overline{\{z\}}$, then $z\in{\rm Im}(f)$.
By the remark above, it is sufficient to prove this
in the case when $y$ is a closed point. For in general,
since $y\in{\rm Im}(f)$, we know that ${\rm Im}(f)\cap \overline{\{y\}}$
is nonempty, so since it is constructible it contains a closed point $y'$.
Then $y' \in \overline{\{y\}} \subseteq \overline{\{z\}}$, 
and hence $z\in{\rm Im}(f)$.

Thus suppose that $z\in Y$, let $Z = \overline{\{z\}}$,
let $y$ be a closed point in ${\rm Im}(f)\cap Z$, and 
assume for a contradiction that $z\notin {\rm Im}(f)$.
Since ${\rm Im}(f)\cap Z$ does not contain $z$, 
this set must be contained in a proper closed subset $F$ of $Z$.
Choose a closed point $t\in Z\setminus F$.
Now it is known that it is possible to join any two points of an
irreducible variety by a curve, see  \cite{M}, p.56.
Thus there is a morphism $C\to Z$,
where $C$ is a curve, and with the image containing $y$ and $t$.
By passing to the normalization if necessary, we may
assume that $C$ is smooth.

Let $c$ be a closed point of $C$ which is sent to $y$.
The complete local ring of $C$ at $c$ is isomorphic to
$\field[[T]]$, and this gives a morphism $h$
\[
{\rm Spec}(\field[[T]]) \longrightarrow C \longrightarrow Z \longrightarrow Y
\]
sending the closed point of ${\rm Spec}(\field[[T]])$ to $y$.
We complete this to a commutative square as in the
statement of the lemma using the morphism
${\rm Spec}(\field) \to X$ corresponding to a closed point
$x\in X$ with $f(x)=y$.
Thus the hypotheses of the lemma give a morphism
$h'$. 

Let $\xi$ be the generic point of ${\rm Spec}(\field[[T]])$. 
Then $h(\xi)$ is the generic point of the closure of the
image of $C$. Thus $\overline{\{h(\xi)\}}$ contains $t$.
On the other hand $h(\xi) = f(h'(\xi)) \in {\rm Im}(f)$,
and $h(\xi)\in Z$, so $h(\xi)\in {\rm Im}(f)\cap Z$.
Thus $\overline{\{h(\xi)\}} \subseteq F$, which does not
contain $t$. This is a contradiction.
\end{proof}

The following observation will also be useful.
The proof follows from the fact that 
$\eta^{-1}(U)$ is closed under generization.

\begin{Lem}
\label{l:opentrick}
Let $\eta: {\rm Spec}(\field[[T]]) \to X$
be a morphism of schemes. 
If $U$ is an open subset of $X$ which contains
the image of the closed point of ${\rm Spec}(\field[[T]])$,
then the image of $\eta$ is contained in $U$. 
\end{Lem}

\section{Proof of Theorem \ref{th:multi}{\rm (iii)}}

Recall that ${\rm mod}_A^d(\field)$ can be identified with
the set of closed points, or the set of $\field$-valued points,
of a scheme ${\rm mod}_A^d$ over $\field$. It is
constructed as follows. 
Let $a_1, \cdots, a_N$ be a set of $\field$-algebra generators of $A$.
For $d \geq 1$, let
\[ 
L = \field[X_{ijl} \mid 1 \leq i,j \leq d, 1 \leq l \leq N] 
\]
be the polynomial ring over $\field$ in $N d^2$ commuting variables.
By $X_l = (X_{ijl})_{1 \leq i,j \leq d}$ we denote the $d \times d$
matrix with $ij$th entry the variable $X_{ijl}$.
Let $\field[A,d]$ be the quotient of $L$ by the ideal
generated by the entries of $f(X_1, \cdots, X_N)$ for all
noncommutative polynomials $f$ over $\field$
in $N$ variables, with $f(a_1, \cdots, a_N) = 0$. 
Define
\[ 
{\rm mod}_A^d = {\rm Spec}(\field[A,d]). 
\]
One can show that this definition is (up to ${\rm GL}_d$-equivariant 
isomorphism) independent of the choice of the generators of $A$,
so ${\rm mod}_A^d$ is called {\it the scheme of} $d$-{\it dimensional}
$A$-{\it modules}.
See \cite{B1}, \cite{Ga} or \cite{P} for the basic properties of schemes of modules.

For a finitely generated commutative $\field$-algebra $R$,
let ${\rm mod}_A^d(R)$ be the set of $R$-valued points of the scheme
${\rm mod}_A^d$.
These are by definition the morphisms
\[ 
{\rm Spec}(R) \longrightarrow {\rm Spec}(\field[A,d]) =  
{\rm mod}_A^d
\]
of schemes over $\field$.
They are in 1-1 correspondence with the $\field$-algebra
homomorphisms
\[ 
\field[A,d] \longrightarrow R. 
\]
Such homomorphisms correspond to $N$-tuples of matrices in $\M_d(R)$
which satisfy the same relations as the generators $a_i$.
Thus they correspond to $\field$-algebra homomorphisms
\[ 
\mathbf{m}:A \longrightarrow \M_d(R). 
\]
Thus the $R$-valued points of ${\rm mod}_A^d$
correspond to $R$-$A$-bimodule structures on $R^d$.

By ${\rm GL}_d$ we denote the general linear 
affine group scheme. It acts on ${\rm mod}_A^d$ by conjugation,
$(\mathbf{g}\cdot\mathbf{m})(a) = \mathbf{g}\mathbf{m}(a)\mathbf{g}^{-1}$
for $\mathbf{g}\in {\rm GL}_d(R)$ and $\mathbf{m}\in {\rm mod}_A^d(R)$.
The orbits correspond to isomorphism classes of $R$-$A$-bimodules 
which are free of rank $d$ over $R$.

Now let $d=d_1+d_2$, let $S$ be a $G$-stable open subset of
${\rm mod}_A^{d_1}(\field)\times {\rm mod}_A^{d_2}(\field)$, and
assume that ${\rm Ext}^1_A(M_2,M_1)=0$ whenever $(M_1,M_2)$ is in $S$.

Let $Z$ be the closed subscheme of ${\rm mod}_A^d$
whose $R$-valued points are the homomorphisms $\mathbf{m}:A\to \M_d(R)$
which take an upper triangular block form with respect to
the decomposition $d=d_1+d_2$. There is a morphism of
schemes over $\field$,
\[
\Delta : Z\longrightarrow {\rm mod}_A^{d_1}\times {\rm mod}_A^{d_2}
\]
which on $R$-valued points sends
an upper triangular homomorphism $\mathbf{m}:A\to \M_d(R)$
to the pair $(\mathbf{m}_1,\mathbf{m}_2)$ consisting
of the two diagonal blocks of $\mathbf{m}$.
Now $S$ can be considered as an open subscheme of 
${\rm mod}_A^{d_1}\times {\rm mod}_A^{d_2}$,
and so there is a corresponding open subscheme
$\Delta^{-1}(S)$ of $Z$.

Let
\[
f : {\rm GL}_d \times \Delta^{-1}(S) \longrightarrow {\rm mod}_A^d
\]
be the morphism obtained from the inclusion of $\Delta^{-1}(S)$ in $Z$,
followed by the conjugation action of ${\rm GL}_d$.
We show that the hypotheses of Lemma \ref{l:schtrick}
are satisfied by the morphism $f$. 
Thus suppose we are given morphisms $p$ and $h$ forming a 
commutative square as in the lemma.

The morphism $h$ is a $\field[[T]]$-valued point of
${\rm mod}_A^d$, so it corresponds to a homomorphism $\mathbf{m}:A\to \M_d(\field[[T]])$.

The morphism $p$ corresponds to a pair $(g_0,m)$ where $g_0\in {\rm GL}_d(\field)$ and
$m$ is a homomorphism $A\to \M_d(\field)$ which has upper triangular
block form, and such that the diagonal blocks satisfy $(m_1,m_2)\in S$.

The commutativity of the square implies that $\mathbf{m}$ is a deformation of $g_0\cdot m$,
so that $g_0^{-1}\cdot \mathbf{m}$ is a deformation of $m$.
Applying Lemma \ref{l:deform} to this deformation gives an 
element $\mathbf{g}\in {\rm GL}_d(\field[[T]])$
such that $\mathbf{m}' = (\mathbf{g} g_0^{-1})\cdot \mathbf{m}$ is a 
homomorphism $A\to \M_d(\field[[T]])$
which has upper triangular block form, and is a deformation of $m$.

Now $\mathbf{m}'$ defines a morphism $\theta:{\rm Spec}(\field[[T]]) \to Z$.
The closed point of ${\rm Spec}(\field[[T]])$ is sent to $m$, and then under $\Delta$
this is sent to $(m_1,m_2)\in S$. Thus Lemma \ref{l:opentrick} implies that
the composition $\theta\Delta$ has image contained in $S$, so that
$\theta$ defines a morphism $\theta':{\rm Spec}(\field[[T]]) \to \Delta^{-1}(S)$.

Now $g_0 \mathbf{g}^{-1}\in {\rm GL}_d(\field[[T]])$ defines a morphism 
${\rm Spec}(\field[[T]]) \to {\rm GL}_d$.
Combining it with $\theta'$ gives a morphism
\[
h' : {\rm Spec}(\field[[T]]) \longrightarrow {\rm GL}_d \times \Delta^{-1}(S).
\]
Since $(g_0 \mathbf{g}^{-1}) \cdot ((\mathbf{g} g_0^{-1})\cdot \mathbf{m}) = 
\mathbf{m}$, we have $h'f = h$.
We also have $\pi^*h' = p$.
Thus the hypotheses of Lemma \ref{l:schtrick} are satisfied,
so the image of $f$, which is ${\mathcal E}(S)$, is open.

\section{Proof of Corollary \ref{co:min}}

Let $C_1$ and $C_2$ be irreducible components of 
${\rm mod}_A^{d_1}(\field)$ and ${\rm mod}_A^{d_2}(\field)$ respectively,
and suppose that ${\rm ext}_A^1(C_2,C_1) = 0$. Let $d=d_1+d_2$.

Let $S$ be the subset of $C_1\times C_2$ corresponding
to the modules $(M_1,M_2)$ with the property that 
$\dim {\rm Ext}_A^1(M_2,M_1) = 0$ and
$\dim {\rm Ext}_A^1(M_1,M_2) = {\rm ext}_A^1(C_1,C_2)$.
This is an irreducible locally closed subset of
${\rm mod}_A^{d_1}(\field)\times {\rm mod}_A^{d_2}(\field)$,
so $\overline{{\mathcal E}(S)}$ is irreducible by
Theorem \ref{th:multi}(ii).

Let $C_i^0$ be the set of elements of $C_i$ which are
contained in no other irreducible component of ${\rm mod}_A^{d_i}(\field)$.
This is an open subset of ${\rm mod}_A^{d_i}(\field)$.
Let $T = S \cap (C_1^0 \times C_2^0)$. This is a nonempty open
subset of ${\rm mod}_A^{d_1}(\field)\times {\rm mod}_A^{d_2}(\field)$,
so by Theorem \ref{th:multi}(i), ${\mathcal E}(T)$ is an
open subset of ${\rm mod}_A^d(\field)$, clearly nonempty.
Thus $\overline{{\mathcal E}(S)}$ contains a nonempty
open subset of ${\rm mod}_A^d(\field)$, so it must be an irreducible
component.

Finally, ${\mathcal E}(C_1\times C_2)$ contains 
$\overline{{\mathcal E}(S)}$ by Theorem \ref{th:multi}(i).

\section{Proof of Corollary \ref{co:filt}}

Using the fact that ${\rm Ext}_A^1(M_i,M_j) = 0$
for all $i \leq j$, an induction shows that any module $N$ in 
${\mathcal F}(M_1, \cdots, M_t)$ 
has a filtration
\[
0 = N_0 \subseteq N_1 \subseteq \dots \subseteq N_t = N
\]
with each $N_j/N_{j-1}$ isomorphic to a direct sum of copies of $M_j$.
Assume by induction on $t$ that ${\mathcal F}_d(M_1, \cdots, M_{t-1})$
is open for all $d$. 
Let $\dim M_t = r$, and for $i\ge 0$ let ${\mathcal O}_i$ be
the orbit in ${\rm mod}_A^{ri}(\field)$ corresponding to the
direct sum of $i$ copies of $M_t$. Then
\[
{\mathcal F}_d(M_1, \cdots, M_t) = \bigcup_{i\ge 0} 
{\mathcal E}({\mathcal O}_i \times {\mathcal F}_{d-ri}(M_1, \cdots, M_{t-1}))
\]
which is open by Theorem \ref{th:multi}(iii), since
${\rm Ext}^1_A(N,M_t)=0$ for all $N$ in ${\mathcal F}(M_1, \cdots, M_{t-1})$.

\section{Proof of Theorem \ref{th:sum}}

Let $C_1, \cdots, C_t$ be irreducible
components of varieties of $A$-modules.
Assuming that 
\[
C = \overline{C_1\oplus\dots\oplus C_t}
\]
is an irreducible component of ${\rm mod}_A^d(\field)$,
we prove that ${\rm ext}_A^1(C_i,C_j) = 0$ for all $i \not= j$.

We say that a $d$-dimensional $A$-module $M$ is 
a {\it minimal degeneration} if the corresponding orbit 
${\mathcal O}(M)$ in ${\rm mod}_A^d(\field)$ is
not contained in the closure of any other orbit.
Note that $\dim {\mathcal O}(M) = d^2 - \dim {\rm End}_A(M)$,
so if an orbit ${\mathcal O}(M)$ is contained
in the closure of another orbit ${\mathcal O}(N)$,
then $\dim {\rm End}_A(M) > \dim {\rm End}_A(N)$.

Recall that the function ${\rm mod}_A^d(\field)\to \mathbb{Z}$, 
sending a module structure to the dimension of the
endomorphism algebra of the corresponding module,
is upper semicontinuous.
Thus if $C$ is an irreducible component of a variety of
modules, then the set $C^{\rm min}$ 
where this function takes its minimal value (on $C$)
is a dense open subset of $C$. 
If $C^0$ denotes the set of elements of $C$
which are not contained in any other irreducible component of 
${\rm mod}_A^d(\field)$, then $C^{0,\rm min} = C^0\cap C^{\rm min}$
is also a a dense open subset of $C$.

Clearly any module in $C^{0,\rm min}$ is a minimal degeneration,
for if ${\mathcal O}(M)$ is contained in the closure of
another orbit ${\mathcal O}(N)$ then the fact that $M$ is in $C^0$
implies that $N$ is in $C$, and then 
$\dim {\rm End}_A(M) > \dim {\rm End}_A(N)$.

Now assume that ${\rm ext}_A^1(C_i,C_j) \not= 0$ for some $i \not= j$.
Thus ${\rm Ext}_A^1(X,Y) \not= 0$ for all modules $X$ in $C_i$ and $Y$ 
in $C_j$.
Taking a non-split short exact sequence 
\[
0 \longrightarrow Y \longrightarrow E \longrightarrow X \longrightarrow 0
\]
we have $E \not\cong X\oplus Y$ by \cite{Ri}, \S 2.3, Lemma 1, so
that $E\oplus Z\not\cong X\oplus Y\oplus Z$ for any module $Z$,
and
${\mathcal O}(X\oplus Y\oplus Z) \subset \overline{{\mathcal O}(E\oplus Z)}$
by \cite{B2}.
It follows that no element of $C_1\oplus\dots\oplus C_t$ can be a 
minimal degeneration. 

But this set must meet $C^{0,\rm min}$. Contradiction.

We now prove the converse direction. By induction
it suffices to prove the case when $t=2$, that is,
if $C_1$ and $C_2$ are irreducible components of 
${\rm mod}_A^{d_1}(\field)$ and ${\rm mod}_A^{d_2}(\field)$, and if 
${\rm ext}_A^1(C_1,C_2) = {\rm ext}_A^1(C_2,C_1) = 0$,
then $\overline{C_1\oplus C_2}$
is an irreducible component of ${\rm mod}_A^d(\field)$, where
$d = d_1+d_2$. 

Let $S$ be the subset appearing in the statement of Corollary \ref{co:min}.
If $(M_1,M_2)$ is in $S$, then ${\rm Ext}^1_A(M_1,M_2)=0$,
so that ${\mathcal E}(S) \subseteq C_1\oplus C_2$.
Now the closure of ${\mathcal E}(S)$ contains an irreducible component
of ${\rm mod}_A^d(\field)$ by Corollary \ref{co:min}.
Thus, since $\overline{C_1\oplus C_2}$ is irreducible,
it must be this irreducible component.


\section{Remarks and examples}\label{reduced}


\subsection{}
There are only very few algebras $A$ such that
all irreducible components of ${\rm mod}_A^d(\field)$
are known for all $d$.
It should be important to construct more explicit 
examples.
If $A$ is of finite representation type, i.e.\ there are only 
finitely many isomorphism classes of indecomposable $A$-modules,
then one can use Auslander-Reiten theory and \cite{Z} to
compute all irreducible components of varieties of
$A$-modules.
One of the known examples of infinite representation type
can be found in \cite{Sch}.

\subsection{}
Theorem \ref{th:sum} can be proved using tangent space methods,
rather than deformation theory, when the scheme ${\rm mod}_A^d$ is 
reduced, or at least generically reduced.
It should be an important task to find methods that
determine when this happens. 
If $A$ is hereditary, i.e.\ the global dimension of $A$ is one,
then ${\rm mod}_A^d$ is reduced. 
One of the few other examples where this is known to be the case
is the variety of complexes, studied by De Concini and Strickland \cite{CS}. 
Namely, let $Q$ be the quiver
\[
m \longrightarrow \dots \longrightarrow 2 \longrightarrow 1 \longrightarrow 0
\]
and let $\Lambda = \field Q/I$, where $I$ is the ideal generated by all 
paths of length 2.
Let $\mathbf{d} = (n_m,\dots,n_1,n_0)$ be a dimension vector,
and define $k_i = \min(n_{i-1},n_i)$. Then ${\rm mod}_\Lambda^{\mathbf{d}}$, 
the scheme of $\Lambda$-modules of dimension vector $\mathbf{d}$, can be 
identified with 
${\rm Spec}(A/{\mathcal E}(k_1,\dots,k_m))$, as in the introduction to 
\cite{CS}, 
and Theorem 1.7 of that paper says this scheme is reduced.

\subsection{}
As mentioned in the introduction, our main results imply that the 
classification of
irreducible components of ${\rm mod}_A^d(\field)$ can
be reduced to the study of irreducible components that
contain dense subsets of indecomposable modules.
So we suggest the following definition.
The {\it component quiver} ${\mathcal C}(A)$ of $A$ has as vertices 
the set of irreducible components $C$ of varieties 
of $A$-modules such that $C$ contains a dense subset of 
indecomposable $A$-modules.
We draw an arrow $C_1 \to C_2$ if 
${\rm ext}_A^1(C_2,C_1) = 0$.

The {\it component graph} ${\mathcal C}_\Gamma(A)$ is the graph
with the same vertices as ${\mathcal C}(A)$
and an edge connecting $C_1$ and $C_2$ if
${\rm ext}_A^1(C_1,C_2) = {\rm ext}_A^1(C_2,C_1) = 0$.
One can easily construct examples where ${\mathcal C}_\Gamma(A)$
contains finite and infinite connected components.

\subsection{}
We give an example to show that ${\mathcal E}(C_1\times C_2)$
in Corollary \ref{co:min} need not be a union of irreducible components.
Let $Q$ be the quiver
\[
1 \xrightarrow{\:\alpha\:} 2 \xrightarrow{\:\beta\:} 
3 \xrightarrow{\:\gamma\:} 4
\]
and let $A = \field Q/I$ where $I$ is generated by the paths $\alpha\beta$
and $\beta\gamma$. The variety ${\rm mod}^d_A(\field)$ breaks into
connected components, one for each decomposition $d=d_1+d_2+d_3+d_4$,
with the component corresponding to the representations of $Q$ of
dimension vector $(d_1,d_2,d_3,d_4)$ which satisfy the relations
$\alpha\beta=0$ and $\beta\gamma=0$.

For the variety ${\rm mod}^1_A(\field)$ there are four 
connected components, each of which is actually a point.
Let $C_1$ be the connected component for the dimension
vector $(1,0,0,0)$.

The variety ${\rm mod}^3_A(\field)$ has a connected component corresponding
to the representations of dimension vector $(0,1,1,1)$. It consists of two 
irreducible 
components, one consisting of the representations with $\beta=0$, the
other corresponding to the representations with $\gamma=0$. Let $C_2$ be
the second of these.

Then ${\mathcal E}(C_1\times C_2)$ consists of the representations of
dimension vector $(1,1,1,1)$ in which $\gamma=0$. It is the union of two
irreducible closed subsets. The subset of representations with $\alpha=0$ is an
irreducible component of ${\rm mod}^4_A(\field)$. The subset of representations
with $\beta=0$ is not an irreducible component, since it is a proper subset
of the set of representations with $\beta=0$ still, but with both 
$\alpha$ and $\gamma$ arbitrary.

\subsection{}
Let $Q$ be the quiver with two vertices $1$ and $2$ and
arrows $\alpha: 1 \rightarrow 1$, $\beta: 1 \rightarrow 1$
and $\gamma: 1 \rightarrow 2$.
Define $A = \field Q/I$ where the ideal $I$ is generated
by $\alpha \gamma$ and $\beta \gamma$.
Then ${\rm mod}_A^1(\field)$ has two connected components,
$C_1$ and $C_2$, corresponding to the representations of dimension
vector $(1,0)$ and $(0,1)$. They are both irreducible, with
${\rm dim}(C_1) = 2$ and ${\rm dim}(C_2) = 0$,
and ${\rm ext}_A^1(C_2,C_1) = 0$.
One can show that 
${\mathcal E}(C_1\times C_2)$ is the union of two 
irreducible components of ${\rm mod}_A^2(\field)$,
of dimensions $3$ and $4$.

\subsection{}
The condition that ${\rm ext}_A^1(C_2,C_1) = 0$ is not necessary in Corollary \ref{co:min}.
Let $Q$ be the quiver with two vertices $1$ and $2$ and
arrows $\alpha: 1 \rightarrow 2$ and $\beta: 2 \rightarrow 1$.
Let $A = \field Q/I$ where $I$ is generated by $\alpha\beta$
and $\beta\alpha$.
Then ${\rm mod}_A^1(\field)$ has two connected components,
$C_1$ and $C_2$, corresponding to the representations of dimension
vector $(1,0)$ and $(0,1)$. They both consist of one point.
One easily checks 
that ${\rm ext}_A^1(C_1,C_2)$ and ${\rm ext}_A^1(C_2,C_1)$ are both
non-zero.
But ${\mathcal E}(C_1 \times C_2)$ and ${\mathcal E}(C_2 \times C_1)$ are 
both irreducible components
of ${\rm mod}_A^2(\field)$ of dimension $3$.

\subsection{}
We would like to add an application to tilted algebras that we worked out 
with C. Gei\ss{} after the paper was submitted.

Let $A=kQ$ be the path algebra of a quiver $Q$.
Recall that the irreducible components of ${\rm mod}_A^d(k)$ are 
indexed by the dimension vectors $\alpha$ for $Q$ with total dimension $d$.
The numbers ${\rm hom}_A(\alpha,\beta)$ and ${\rm ext}^1_A(\alpha,\beta)$
may be computed using the theory of general representations of quivers \cite{S},
and they are related by the formula
${\rm hom}_A(\alpha,\beta)-{\rm ext}^1_A(\alpha,\beta) = \langle\alpha,\beta\rangle$
where $\langle -,- \rangle$ is the Ringel form for $Q$.

Let $T_A$ be a tilting module, let $t$ be its dimension vector,
and let $B = {\rm End}_A(T)^{op}$ 
be the corresponding tilted algebra, see \cite{HR}.
(The opposite algebra is needed, with our conventions, 
to ensure that $T$ is naturally a $B$-$A$-bimodule.) 
We show that the irreducible components of ${\rm mod}_B^d(k)$ 
are in 1-1 correspondence with the pairs of dimension vectors 
$(\alpha, \beta)$ for $Q$ with 
$\langle t,\alpha-\beta\rangle = d$, 
${\rm ext}_A^1(t,\alpha)=0$, 
${\rm hom}_A(t,\beta)=0$ and ${\rm hom}_A(\beta,\alpha)=0$. 

The Brenner-Butler Theorem \cite[Theorem 2.1]{HR} gives equivalences 
$F = {\rm Hom}_A(T,-):\mathcal{T}\to\mathcal{Y}$ and 
$F' = {\rm Ext}_A^1(T,-):\mathcal{F}\to\mathcal{X}$ involving the module classes
\begin{align*}
\mathcal{F} &= \{ \text{$M$ an $A$-module} \mid {\rm Hom}_A(T,M) = 0\},
\\
\mathcal{T} &= \{ \text{$M$ an $A$-module} \mid {\rm Ext}^1_A(T,M) = 0\},
\\
\mathcal{X} &= \{ \text{$N$ a $B$-module} \mid N\otimes_B T = 0\} = \{ N \mid {\rm Hom}_B(N,DT)=0 \},
\\
\mathcal{Y} &= \{ \text{$N$ a $B$-module} \mid {\rm Tor}_1^B(N,T) = 0\} = \{ N \mid {\rm Ext}^1_B(N,DT)=0 \}.
\end{align*}
By upper semicontinuity, $\mathcal{F}$ and $\mathcal{T}$ 
define open subsets $\mathcal{F}(d)$ and $\mathcal{T}(d)$ of 
${\rm mod}_A^d(k)$, and $\mathcal{X}$ and $\mathcal{Y}$ define open 
subsets $\mathcal{X}(d)$ and $\mathcal{Y}(d)$ of ${\rm mod}_B^d(k)$.

The irreducible components of $\mathcal{T}(d)$ are the ones for 
${\rm mod}_A^d(k)$ containing a module in $\mathcal{T}$, and since $T$ corresponds to
an open orbit, they are the ones with ${\rm ext}_A^1(t,\alpha)=0$. 
By a result of Bongartz \cite[Theorem 3]{Bmin}, 
the modules $F(M)$ with $M$ in 
such an irreducible component form an irreducible component of $\mathcal{Y}(d_1)$, 
where $d_1=\langle t,\alpha \rangle$. 
Varying $\alpha$, this construction gives all irreducible components of $\mathcal{Y}(d_1)$.

Dually, the dimension vectors $\beta$ with ${\rm hom}_A(t,\beta)=0$ index the irreducible 
components of $\mathcal{X}(d_2)$, where $d_2 = - \langle t,\beta \rangle$.
(This requires the dual of Bongartz's Theorem, but this follows automatically 
on exchanging the roles of $A$ and $B$, and using the dual of $T$.)

Now any $B$-module is a direct sum of modules in $\mathcal{X}$ and $\mathcal{Y}$ by
\cite[Theorem 6.3]{HR}, so any irreducible component of 
${\rm mod}_B^d(k)$ is of the form 
$\overline{C_1\oplus C_2}$ where $C_1$ is the closure in ${\rm mod}_B^{d_1}(k)$ 
of an irreducible component of $\mathcal{Y}(d_1)$, say corresponding to $\alpha$, 
and $C_2$ is the closure in ${\rm mod}_B^{d_2}(k)$ 
of an irreducible component of $\mathcal{X}(d_2)$, say corresponding to $\beta$. 
Moreover, the summands of the $B$-module in $\mathcal{X}$ and $\mathcal{Y}$ 
are unique up to isomorphism, and it follows that different choices of 
$C_1$ and $C_2$ lead to different irreducible components.

Now the condition for $\overline{C_1\oplus C_2}$ to be an irreducible component is that ext 
vanishes generically in both directions. In one direction this is automatic, for 
${\rm Ext}_B^1(\mathcal{Y},\mathcal{X}) = 0$ by \cite[4.1.6(d)]{Ri}. In the other direction we have
\[
{\rm Ext}_B^1(F'(M) ,F(N)) \cong {\rm Hom}_A(M,N)
\]
for $M\in\mathcal{F}$ and $N\in\mathcal{T}$. (This is obvious from the
interpretation of tilting theory as an equivalence of derived categories.)
The condition for this to vanish generically is that ${\rm hom}_A(\beta,\alpha)=0$.

\end{document}